\documentclass[runningheads]{llncs}
\usepackage{graphicx}
\usepackage{url}
\usepackage{amsmath}
\usepackage{amssymb}
\usepackage{booktabs}
\usepackage{tabularx}
\usepackage{setspace}
\usepackage{lscape}
\usepackage{float}
\usepackage{algorithmicx}
\usepackage{algorithm}
\usepackage{algcompatible}
\usepackage{algpseudocode}
\usepackage{listings}
\usepackage{breqn}
\usepackage{color}
\usepackage{comment}

\def \dj{d\kern-0.3em\char"16\kern-0.3em}
\def \Dj{\mbox{\raise0.3ex\hbox{-}\kern-0.4em D}}

\begin{document}

\title{Metaheuristics for finding threshold graphs with maximum spectral radius}
\titlerunning{MH for Max Sp Radius Threshold Graphs}
\author{Luka Radanovi\'c\inst{1}, Abdelkadir Fellague\inst{2}, Dragutin Ostoji\'c\inst{3}, 
Dragan Stevanovi\'c\inst{4}\thanks{DS is on a sabbatical leave from the Mathematical Institute of the Serbian Academy of Sciences and Arts}, Tatjana  Davidovi\'{c}\inst{5}}
\authorrunning{Radanovi\'c et al.}
%

\institute{
Faculty of Mathematics, University of Belgrade, Serbia,
\email{lukaradanovich@gmail.com}
\and
SIMPA, Department of Computer Science, University of Sciences and Technology of Oran, Algeria,
\email{ abdelkadir.fellague@univ-usto.dz}
\and
Faculty of Science and Mathematics, University of Kragujevac, Serbia,
\email{dragutin.ostojic@pmf.kg.ac.rs}
\and
Abdullah Al Salem University, Khaldiya, Kuwait,
\email{dragan.stevanovic@aasu.edu.kw}
\and
Mathematical Institute of the Serbian Academy of Sciences and Arts, Belgrade, Serbia,
\email{tanjad@mi.sanu.ac.rs}
}

\maketitle              
\date{}
\begin{abstract}
We consider the problem of characterizing graphs with the maximum spectral radius among the connected graphs with given numbers of vertices and edges.
It is well-known that the candidates for extremal graphs are threshold graphs, but only a few partial theoretical results have been obtained so far.
Therefore, we approach to this problem from a novel perspective that involves incomplete enumeration of different threshold graphs with a given characteristic. Our methodology defines the considered problem as an optimization task and utilizes two metaheuristic methods, Variable Neighborhood Search (VNS), which relies on iterative improvements of a single current best solution and Bee Colony Optimization (BCO), a population-based metaheuristic from the Swarm Intelligence (SI) class. We use compact solution representation and several auxiliary data structures that should enable efficient search of the solution space. In addition, we define several types of transformations that preserve the feasibility of the resulting solution. The proposed methods are compared on the graphs with a moderate number of vertices. Preliminary results are in favor of the VNS approach, however, we believe that both methods could be improved.

\medskip

\noindent {\bf Keywords:} Spectral graph theory; extremal graphs; spectral radius; adjacency matrix; metaheuristics.

\end{abstract}

\section{Introduction}

Graphs are mathematical objects defined as 2-tuples  $G = (V,E)$ \cite{Cve90}, where $V = \{v_1, v_2, \ldots, v_n\}$, represents the set of \emph{vertices}
$v_i$, while $E\subseteq V\times V$ denotes the connections (relations) between the pairs of vertices and is called the set of \emph{edges}. If there is a connection (edge) between vertices $v_i$ and $v_j$, we say that $\{v_i, v_j\}\in E$ and that vertices $v_i$ and $v_j$ are \emph{adjacent}. Graphs are used to model numerous problems in science, engineering, industry, etc. Usually, $V$ is finite set, however, the infinite cases are also studied in the literature starting with \cite{Nas67}. In this paper, we consider only finite and undirect graphs with a pre-specified number of vertices and edges. 


The simplest graph representation is by the \emph{Adjacency matrix} $A$ with elements 0 or 1 defined as follows:
$$
a_{ij} = \left\{\begin{array}{ll}
1, & \mbox{if $\{v_i, v_j\}\in E$};\\
0,& \mbox{otherwise}.
\end{array}\right.
$$
If graph is undirected, $A$ is symmetric, i.e., $a_{ij}=a_{ji}. $The \emph{degree} of vertex $v_i$ (denoted by $d_{i}$) in graph $G$ represents the number of vertices adjacent to $v_i$, i.e., the number of edges having $v_i$ as an end-vertex and it is calculated as $d_{i} = \sum_{j=1}^n a_{ij}$. \emph{Eigenvalues} $\lambda_i$, $i=1,2,\ldots, n$ for the graph $G$ are actually the eigenvalues of matrix $A$, i.e., the roots of its characteristic polynomial $P_G(x)=det(x I - A)$. As the adjacency matrix $A$ is symmetric its eigenvalues are real numbers. The set of all eigenvalues of graph $G$ is called \emph{spectrum}. It can contain negative, positive values and zeros, with some repeated values. It is usual to represent the spectrum as a non-increasing array of values $\lambda_1 \geq \lambda_2 \geq \cdots \geq \lambda_n$. Then, the largest eigenvalue $\lambda_1$ of graph $G$ is called \emph{index}. An array $x$ such that $Ax = \lambda x$ is known as \emph{eigenvector} (corresponding to the eigenvalue $\lambda$) of graph $G$, and it actually represents the eigenvector of matrix $A$.

\emph{Spectral graph theory} (SGT) \cite{CDS95,Spi12} studies graphs based on their adjacency matrix, more precisely, based on the \emph{eigenvalues} and \emph{eigenvectors} of this matrix. In recent literature, some other matrices associated with graphs are defined and analyzed, such as Laplacian matrix and signless Laplacian matrix (\cite{CDS95}, section 1.3). However, they will not be considered in this paper. SGT has important applications in various fields of computer science \cite{CS11}, some of them including finding extremal graphs with respect to a given invariant or the combination of several invariants  \cite{ABFCHHLM06,Bol04,CH00,Erd64} as graphs represent natural models for various types of objects. 

The problem considered in this study is to characterize graphs with the maximum spectral radius (the largest eigenvalue) among connected graphs with given numbers of vertices and edges, which is open for more than 35 years. This problem of characterizing (not necessarily connected) graphs with maximum spectral radius having a given numbers of vertices ($n$) and edges ($m$) was posed initially by Brualdi and Hoffman in 1976 \cite{BFLS78}, (p.~438). The first theoretical results were related to the disconnected graphs, while the problem of characterizing connected extremal graphs is still unresolved in the general case. Brualdi and Solheid \cite{BS86} showed that the adjacency matrix of a connected extremal graph must have a stepwise form, in the sense that its vertices can be ordered in such a way that $A_{ij}=1$ (with $i<j$) implies $A_{hk}=1$ for all $h\leq i$, $k\leq j$ and $h<k$. Simi\'c, Li Marzi and Belardo \cite{SLB04} offered an alternative reasoning showing that a connected extremal graph cannot have either the path $P_4$, the cycle $C_4$ or the pair of independent edges $2K_2$ as an induced subgraph. Namely, this implies that a connected extremal graph is a \emph{threshold graph}.

Threshold graphs can be described iteratively as it is proposed in \cite{MP95}. We start with a single vertex and, in each step, add a new vertex that is either isolated or adjacent to all already included vertices. This process of sequentially building a threshold graph may be written in a more formal way as:
\begin{equation}
G_{p_1} = K_{p_1}
\label{first} 
\end{equation}
\begin{equation}
G_{p_1,p_2,\ldots,p_k}=\overline{G_{p_1,p_2,\ldots,p_{k-1}}} \vee K_{p_k}
\label{k-th} 
\end{equation}
where $p_1, p_2,\ldots, p_k$, are positive integers, $\overline{G}$ denotes the complement of $G$ and $\vee $ denotes the join of two graphs.
This notation compresses successive additions of $p_1$ vertices of one type (each isolated or each adjacent to all previous vertices), $p_2$ vertices of the opposite type, $p_3$ vertices of the first type, etc. Here the complement changes the types of previous vertices, while the join ensures that the $p_k$ vertices added at the last step are adjacent to all previous vertices.

\begin{figure}[!htbp]
  \begin{center}
$$
A=\left[\begin{array}{cccccccc}
0& 0& 0& 1& 0& 1& 0& 1\\
0& 0& 0& 1& 0& 1& 0& 1\\
0& 0& 0& 1& 0& 1& 0& 1\\
1& 1& 1& 0& 0& 1& 0& 1\\
0& 0& 0& 0& 0& 1& 0& 1\\
1& 1& 1& 1& 1& 0& 0& 1\\
0& 0& 0& 0& 0& 0& 0& 1\\
1& 1& 1& 1& 1& 1& 1& 0
        \end{array}
\right] \;\;\;\;\; n=8,\; m=15,\; \lambda_1 = 4.37
$$
\hspace{2cm}\includegraphics[width=3cm]{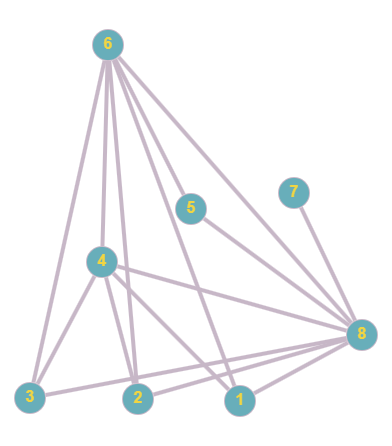}
    \caption{An example of threshold graph}
    \label{fig:exampl}
  \end{center}
\end{figure}

An example of threshold graph with $n=8$ vertices and $m=15$ edges is presented in Fig.~\ref{fig:exampl}.

Although we know that only threshold graphs are valid candidates for the spectral radius maximizer, the theoretical results in a general case are still missing. Exhaustive enumeration of all threshold graphs with a given $n$ and $m$ is an NP-hard problem, especially for the medium number of edges, i.e., for $m$ close to $n(n-1)/4$. Therefore, the application of some incomplete search methods is more than welcome. One of the possible approaches is to use some general purpose software developed for generating conjectures in graph theory, such as GRAPH, 1984 (https://www.mi.sanu.ac.rs/novi\_sajt/research/projects/GRAPH.zip), the three versions of AutoGraphiX, 1997, 2009, 2015 (https://www.autographix.ca), newGraph 2004 (https://www.mi.sanu.ac.rs/newgraph/), or maybe PHOEG 2008 (https://phoeg.umons.ac.be/phoeg), to mention a few. Another popular approach is to develop metaheuristic, tailored for each particular optimization problem. We illustrate the application of Variable Neighborhood Search (VNS) and Bee Colony Optimization (BCO)  metaheuristics. These two methods are representative of two distinct classes of algorithms, mathematical-based singe-solution and nature-inspired population-based, developed by Serbian scientists.

The remainder of this paper is organized as follows. Section~\ref{lit} provides a brief overview of the relevant literature, specifying complexity of the considered problem and some special cases in which the solution is known and there are theoretical proofs of their correctness. 
Section~\ref{MH} contains an overview of developed metaheuristic methods. In Section~\ref{VNS}, a short description of Variable Neighborhood Search (VNS) is provided, while Section~\ref{BCOi} contains the description of Bee Colony Optimization (BCO) algorithm. Details about the implementation of Geeneral VNS (GVNS) and improvement-based BCO (BCOi) for finding threshold graphs with maximum radius are presented in Section~\ref{impl}. The results obtained by applying the implemented methods to some known graphs from the literature, are presented in Section~\ref{exper}. Concluding remarks and guidelines for future work are given in Section~\ref{concl}.

\section{Literature review}
\label{lit}

The theoretical results related to the identification of spectral radius maximizer threshold graphs can be summarized as follows:
\begin{itemize}
\item	for $m=n-1$, the extremal connected graph is star $K_{(n-1,1)}$, as it is proven in \cite{LP73}. To follow the notation from \eqref{k-th}, we keep the $G$-notation, i.e., $K_{(n-1,1)}=G_{(n-1,1)}$
\item	for $m=n+d$ with $d\in\{0,2\}$, the extremal connected graph is $G_{(2+d/2,n-3-d/2,1)}$ \cite{BS86}.
\item	for $m=n+d$ with $d=1$ or fixed $d\geq 6$ and sufficiently large $n$, the extremal connected graph is $G_{(d+1,1,n-3-d,1)}$ \cite{BS86,CR88}.
\item	for $m=n+\binom{d-1}{2}-1$ with fixed $d\geq 5$, the extremal connected graph is either $G_{(d-1,n-d,1)}$ or $G_{(\binom{d-1}{2},1,n-2-\binom{d-1}{2},1)}$ \cite{Bell91}.
\item	for $m=n+\binom{d-1}{2}-2$ with $2n\leq m<\binom{n}{2}-1$, the extremal connected graph is $G_{(2,d-2,n-1-d,1)}$ \cite{ORV02}.
\end{itemize}

The relevant results related to spectral radius of some additional particular classes of graphs are summarized in the book \cite{Ste15}. 

While threshold graphs and their various spectral properties do attract considerable attention among researchers (see, for example, results on computing characteristic polynomial of threshold graphs in \cite{JTT14,LMT19,ADDS20}), we must mention here that no relevant papers on the topic considered here were published in the last ten years, with the last noteworthy result being a few bounds on the spectral radius of threshold graphs in \cite{SBMT10}. This very much signifies the need for a thorough change of methodology in treating this problem.

\section{Metaheuristic methods overview}
\label{MH}

This section contains the brief description of the Variable Neighborhood Search (VNS) and Bee Colony Optimization (BCO) metaheuristic. In particular, we focus on General VNS and improvement-based variant of BCO (BCOi), as they are the most suitable for the application to the considered problem.

\subsection{Variable neighborhood search}
\label{VNS}

Variable Neighborhood Search (VNS) is a trajectory-based metaheuristic method proposed in \cite{MH97}. It uses distances between solutions and employs one or more neighborhood structures to efficiently search the solution space of a considered optimization problem. VNS uses some problem-specific Local Search (LS) procedure in the exploitation phase and changing distances between solutions to ensure the exploration of solution space. The role of exploration (diversification, perturbation) phase is to enable escaping from local optima traps. VNS is widely used optimization tool with many variants and successful applications \cite{HMBM19}, and we use its general variant (GVNS) to search for threshold graphs with maximum radius. The pseudo-code of GVNS is presented by Alg.~\ref{GVNS-pscode}.

The main steps of most VNS variants are: Shaking, Local Improvement, and Neighborhood Change (see Alg.~\ref{GVNS-pscode}). The role of Shaking step is to ensure the diversification of the search. It performs a random perturbation of $x_{best}$ in the given neighborhood and provides a starting solution $x'$ to the next step. Local Improvement phase tries to improve $x'$ by visiting its neighbors with respect to the selected neighborhood(s). After this phase is completed, VNS performs Neighborhood Change step in which it examines the quality of the obtained local optimum $x''$. If it is better than $x_{best}$, the search is concentrated around it (the global best solution $x_{best}$ and the neighborhood index $k$ are updated properly). Otherwise, only $k$ is changed. The three main steps are repeated until a pre-specified stopping criterion is satisfied \cite{HMBM19}.

The main parameter of VNS is $k_{max}$, the maximum number of neighborhoods for Shaking. Actually, the current value of $k$ represents the distance between $x_{best}$ and $x'$ obtained within the Shaking phase. VNS is known as the First Improvement (FI) search strategy because the search is always concentrated around $x_{best}$: as soon as this solution is improved, $k$ is reset to 1.

\begin{algorithm}[!htb]
\caption{Pseudo-code for GVNS method}
\label{GVNS-pscode}
{\footnotesize
\begin{algorithmic}
\Procedure{GVNS}{Problem input data, $k_{max}$, STOP}
   \State {$x_{best} \gets \Call{InitSolution}$}
   \Repeat
      \State ${k \gets 1}$
      \Repeat
           \State {$x' \gets \Call{RandomSolution}{x_{best},{\cal N}_k}$} \Comment{Shaking}
           \State {$x'' \gets \Call{VND}{x'}$} \Comment{Local Improvement}
           \If {$(f(x'') < f(x_{best}))$} \Comment{Neighborhood Change}
               \State {$x_{best} \gets x''$}
               \State {$k \gets 1$}
           \Else
               \State {$k \gets k+1$}
           \EndIf
           \State {$Terminate \gets \Call{StoppingCriterion}{STOP}$}
      \Until {$(k > k_{max} \vee Terminate)$}
   \Until {(Terminate)}
\State $\Call{Return}{x_{best},f(x_{best})}$
\EndProcedure
\end{algorithmic}
}
\end{algorithm}

Contrary to the Basic VNS (BVNS) that employs a single type of neighborhood within the Local Improvement phase, GVNS explores a set of neighborhoods in a deterministic manner defined by the Variable Neighborhood Descent (VND) procedure~\cite{HMTH17}. This procedure systematically searches the specified set of neighborhoods in a predetermined order according to the first improvement strategy. This means that at each improvement of the current solution, the search is re-directed to the first neighborhood. Only in the case when no improvement can be found in the searched neighborhood, the next one will be explored. VND is a deterministic search procedure that completes when all neighborhoods are explored without the improving the current solution. Alg.~\ref{VND-pscode} illustrates the execution of VND procedure.

\begin{algorithm}[H]
\caption{Variable Neighborhood Descent} \label{VND-pscode}
{\footnotesize
\begin{algorithmic}
\Procedure{VND}{$S_x$, $l_{max}$}
\State {$l\gets 1$}
\While{$(l \leq l_{max})$}
    \State {$S_x^{'}\gets \Call{LocalSearch}{S_x, {\cal N}_l}$}
    \State {$l\gets l + 1$}
    \If{$(S_x^{'}$ \text{better than} $S_x$)}
        \State{$S_x\gets S_x^{'}$}
        \State {$l \gets 1$}
    \EndIf
\EndWhile
\State \Call{Return}{$S_x$}
\EndProcedure
\end{algorithmic}
}
\end{algorithm}


\subsection{Bee Colony Optimization}
\label{BCOi}

Bee Colony Optimization (BCO) is a population-based metaheuristic that mimics the foraging process of honeybees in nature \cite{DTS15}. The population consists of artificial bees, each responsible for one solution of the considered problem. During the execution of BCO, artificial bees build (in the constructive BCO variant, BCOc) or transform (in the improvement-based BCOi) their solutions in order to find the best possible with respect to the given objective. The BCO algorithm runs in iterations until a stopping condition is met and the best found solution (the so called global best) is reported as the final one.


\begin{algorithm}[!htb]
\caption{Pseudo-code of the BCO algorithm}
\label{BCO-pscode}
{\footnotesize
\begin{algorithmic}
\Procedure{BCO}{Problem input data, $B,NC,STOP$}
   \State $Terminate \gets 0$
   \While {$!Terminate$}
      \For{$b\gets 1, B$} \Comment{// Determine initial population}
        \State $Solution(b) \gets \Call{GenerateSolution}$ \Comment{// Initial Forward pass}
      \EndFor
      \State $\Call{Update}{x_{best}}$
      \For{$u\gets 1, NC$}
         \State $\Call{Normalization}$ \Comment{// Backward pass}
         \State $U \gets \Call{Loyalty}$
         \For{$b\gets 1, U$}
               \State $\Call{Recruitment}{Solution(b)}$
         \EndFor
         \For{$b\gets 1, B$} \Comment{// Forward pass}
              \State ${\Call{Transform}{Solution(b)}}$
         \EndFor
         \State $\Call{Update}{x_{best}}$
      \EndFor
      \State $Terminate \gets \Call{StoppingCriterion}{STOP}$
   \EndWhile
   \State \Call{Return}{$x_{best},f(x_{best})$}
\EndProcedure
\end{algorithmic}
}
\end{algorithm}

Each BCO iteration contains several execution steps divided into two alternating phases: \emph{forward pass} and \emph{backward pass} (see Alg.~\ref{BCO-pscode}). Within forward passes, all bees explore the search space by applying a predefined number of moves and obtain new population of solutions. Moves are related to building or transforming solutions, depending on the used BCO variant and they explore \emph{a priori} knowledge about the considered problem. When a new population is obtained, the second phase (backward pass) is executed, where the information about the quality of solutions is exchanged between bees. The solution's quality is defined by the corresponding value of the objective function. The next step in backward pass is to select a subset of promising solutions to be further explored by applying \emph{loyalty decision} and \emph{recruitment} steps. Depending on the relative quality of its current solution with respect to the best solution in the current population, each bee decides with a certain probability should it stay \emph{loyal} to that solution and become a \emph{recruiter} that advertises its solution by simulating waggle dance of honeybees \cite{DTS15}. Obviously, bees with better solutions should have more chances to keep their solutions.

The probability that $b$-th bee (at the beginning of the new forward 
pass) is loyal to its previously generated partial/complete solution 
can be defined in a number of ways \cite{J-K16}, the simplest of them 
exploring the normalized value of the objective function:
\begin{equation}
p_b = O_b, \;\;\; b=1,2,\ldots,B
\label{p_loy}
\end{equation}
where:\\
 $O_b$ - denotes the normalized value for the objective function
of partial/complete solution created by the $b$-th bee;\\
 $u$ - counts the forward passes (taking values 1, 2, \ldots, $NC$).
In the case of maximization problem, normalized value of the objective function $C_b$ is
calculated as follows:
 \begin{equation}
 O_b = \frac{C_b - C_{\min}}{{\displaystyle
 C_{\max} - C_{\min}}},\;\;\;b=1,2,\ldots,B
\label{bees_norm2}
 \end{equation}
where $C_{min}$ and $C_{max}$ are the objective function values related 
to worst (minimal) and best (maximal) solutions, respectively, among 
the all bees engaged in the current forward pass. From equation
(\ref{bees_norm2}), it could be seen that if $b$-th bee
partial/complete solution is closer to minimal value of all
obtained solutions ($C_{min}$), i.e., it corresponds to a lower quality 
solution, than its normalized value, $O_b$, will be smaller. On the other 
hand, if the value of the partial/complete solution $C_b$ is larger, 
then its normalized value $O_b$ will also be larger.

\begin{figure}[!htbp]
  \begin{center}
\includegraphics[width=6cm]{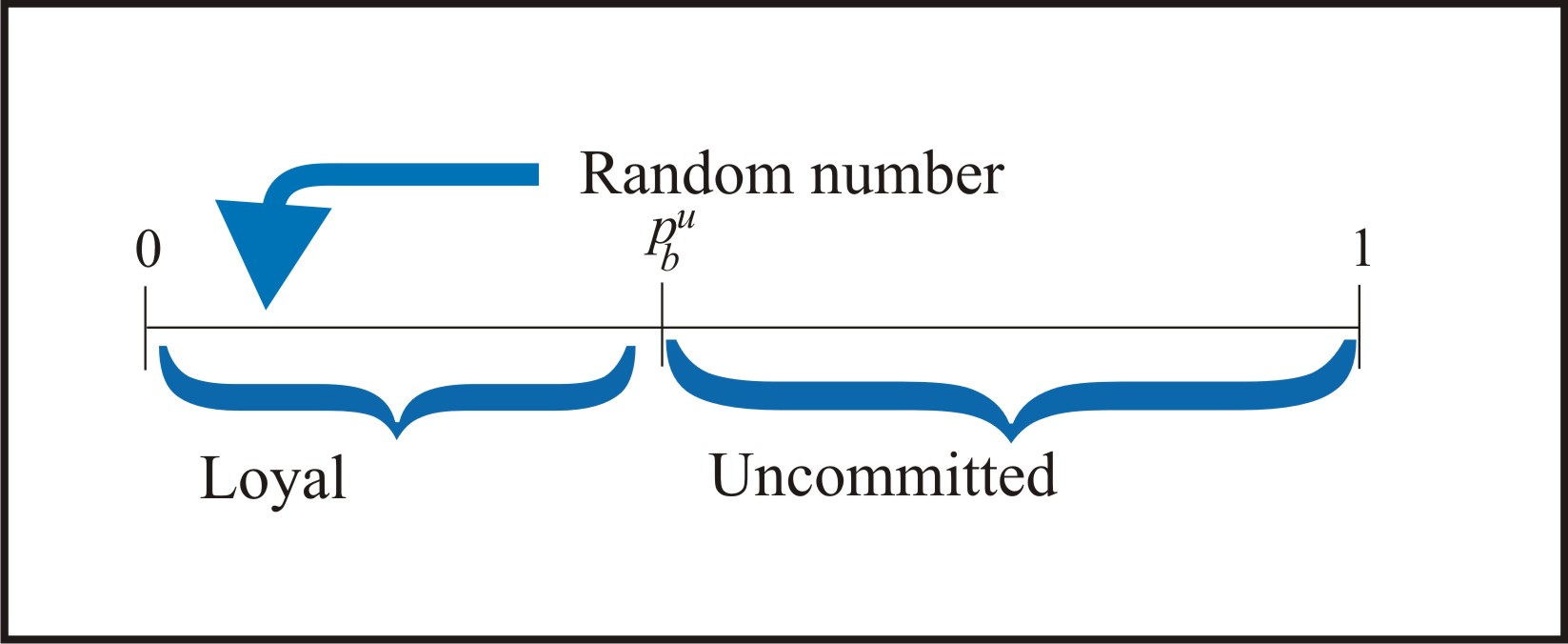}
    \caption{Deciding on loyalty for bee $b$}
    \label{LoyDec}
  \end{center}
\end{figure}

Equation (\ref{p_loy}) and a random number generator are used for
each artificial bee to decide whether it should stay loyal (and
continue exploring its own solution) or to become an
\emph{uncommitted follower}. If the generated random number
from $[0,1]$ interval is smaller than the calculated probability
$p_b$, then artificial bee $b$ stays loyal to its own solution.
Otherwise, the bee becomes uncommitted. The decision on loyalty 
for bee $b$ is illustrated in Fig.~\ref{LoyDec}. 

The uncommitted followers abandon their current solutions and 
have to select one of the solutions held by recruiters that are 
"advertised" for additional exploration. This selection is taken 
with a probability, such that better advertised solutions have 
greater opportunities to be chosen for further exploration. 

\begin{figure}[!htbp]
  \begin{center}
\includegraphics[width=6cm]{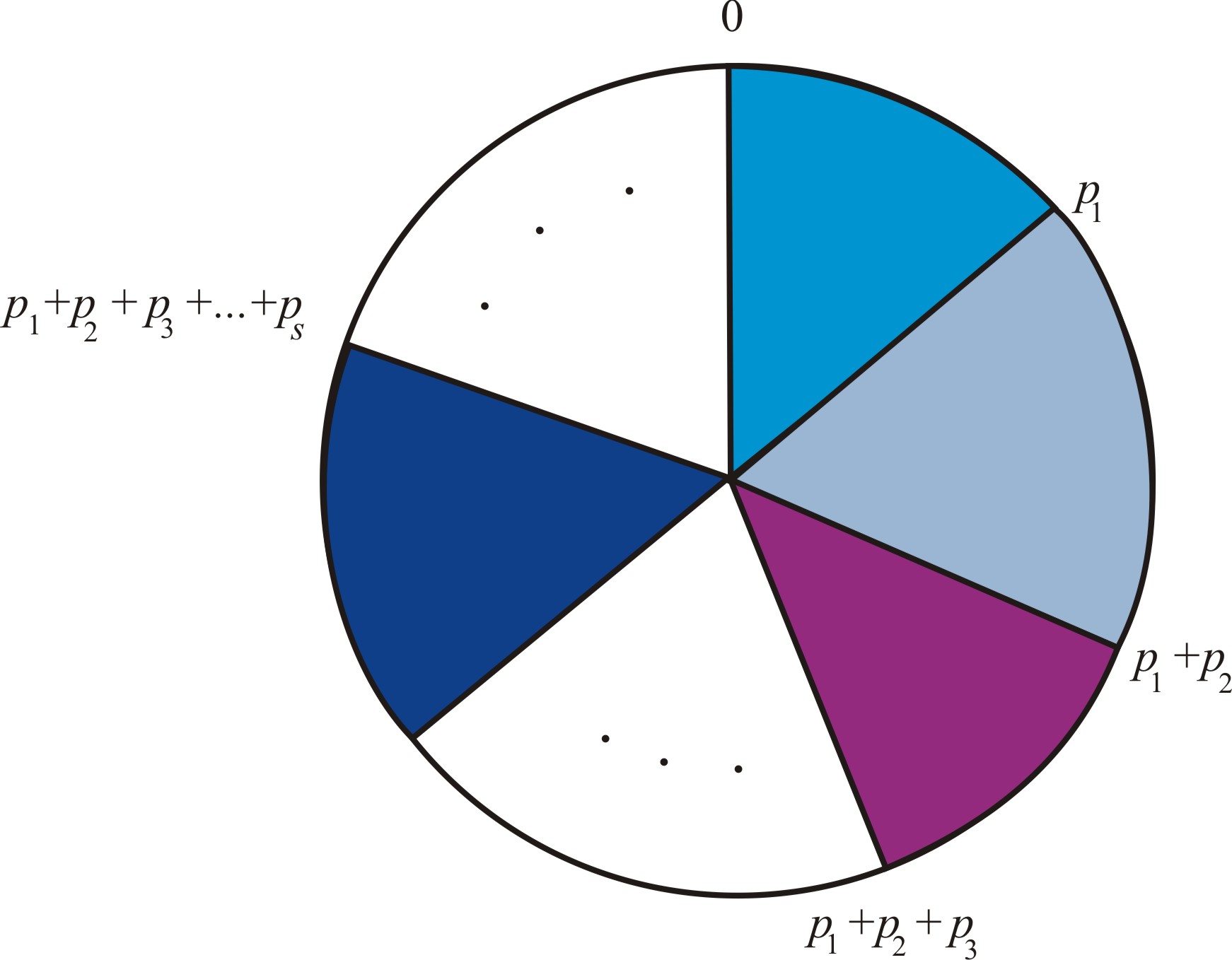}
    \caption{Recruitment of bee $b$ by the roulette wheel}
    \label{roul}
  \end{center}
\end{figure}

For each uncommitted bee, it is decided which recruiter it will
follow, taking into account the quality of all advertised
solutions. The probability that the solution generated by $r$-th recruiter
would be chosen by any uncommitted bee equals:
 \begin{equation}
 p_r = \frac{O_r}{{\displaystyle
 \sum_{k=1}^{R}O_k}},\;\;\;r=1,2,\ldots,R
\label{recr_prob}
 \end{equation}
where $O_k$ represents the normalized value for the objective
function of the $k$-th advertised solution and $R$ denotes the
number of recruiters. Using equation (\ref{recr_prob}) and a
random number generator, each uncommitted follower joins one
recruiter through the various selection mechanisms, the most 
commonly used being roulette wheel (Fig.~\ref{roul}).

The roulette wheel is a well-known model of choice. The main
inspiration for its development came from a game-gambling
roulette. Any solution can be chosen, and the probability of its
selection (the size of a particular slot on the roulette wheel)
depends on the quality of that solution, i.e., on the value of the
objective function.

In practice, the size of the slot on the roulette wheel associated
to each solution is 2determined by the ratio between the corresponding
normalized objective function value and the sum of normalized
objective function values for all advertised solutions. On the one
hand, a solution with better objective function value has a higher
chance to be selected. On the other hand, due to the influence of 
randomness, there is still a possibility that it will be eliminated 
from the continuation of search process.

In the basic 
variant of BCO there are only two parameters:\\
$\bullet$ $B$ -- the number of bees involved in the search (population size) and\\
$\bullet$ $NC$ -- the number of forward/backward passes in a single BCO iteration.

\section{GVNS and BCOi for maximizing index of threshold graphs} 
\label{impl}

To efficiently implement any metaheuristic method, we need to define the solution representation, auxiliary data structures, generation of initial solution and solution transformation rules, i.e., neighborhoods of each particular feasible solution. 

\subsubsection{Solution representation}

A threshold graph on $n$ vertices and $m$ edges may be represented by a binary sequence $\mathbf{R}=\{r_1,r_2,\ldots,r_n\}$, in which $0$ at the position $i$ means that node $i$ is an isolated vertex with respect to all previous vertices, i.e., vertices with indices smaller than $i$, (type 0 vertex), while 1 means that node $i$ is adjacent to all previous vertices (type 1 vertex). The bit value for $r_1$ may be chosen arbitrarily and we set it to 1. As we are searching over the connected threshold graphs, the last vertex should be connected to all the remaining vertices, and therefore, $r_n=1$. In addition, the following equality
$$\sum_{i=1}^{n}(i-1)\cdot r_i = m$$
should hold to ensure that the number of edges equals $m$.

\subsubsection{Initial solution generation}

Initial solution could be generated in a greedy manner. The first element of the solution sequence is initialized to 1, while all other elements are set to 0. In the next step, starting with the last vertex, we add edges from that vertex to all previous if possible. In the case that it is possible to add these edges, the corresponding element in the sequence $\mathbf{R}$ becomes 1 and the number of available edges is decreased by the number of used edges. Otherwise, the value sequence element remains 0. In both cases, we move to the previous vertex. When this loop is completed, we end up with a feasible solution. The pseudo-code of this procedure is presented in Algorithm~\ref{INIT-pscode}.

\begin{algorithm}[!htb]
\caption{Pseudo-code of the INIT\_SOL algorithm}
\label{INIT-pscode}
{\footnotesize
\begin{algorithmic}
\Procedure{INIT\_SOL}{$n,m,R$}
\State $R[1] \gets 1$
\State $R[n] \gets 1$
\For{$i\gets 2, n-1$} 
   \State $R[i] \gets 0$
\EndFor
\State $e = m-n+1$
\State $i\gets n$
\While {$(i > 1) \wedge (e > 0)$}
    \If {$i - 1 <= e$}
        \State {$R[i] \gets 1$}
        \State {$e \gets e - i+1$}
    \EndIf
    \State $i\gets i-1$
\EndWhile
\EndProcedure
\end{algorithmic}
}
\end{algorithm}

However, as BCOi is a population-based method, we need a set of feasible solutions at the beginning of each iteration. Obtaining different feasible solutions from the initial, obtained by the described greedy procedure is easy: it is just necessary to perform a random number of stochastic transformations described in the text to follow. 

However, our preliminary experiments showed that this procedure is not appropriate for having a diverse population of solutions that could be efficiently explored by the BCOi algorithm. Namely, the number of possible transformation is small, therefore, all solutions are similar. The reason for such situation is the main characteristics of this greedy initial solution is that elements of the sequence that have value 1 are usually concentrated together and close to the sequence's tail. Therefore, we developed an alternative procedure for generating an initial solution described by Algorithm~\ref{INIT2-pscode}.

\begin{algorithm}[!htb]
\caption{Pseudo-code of the INIT2\_SOL algorithm}
\label{INIT2-pscode}
{\footnotesize
\begin{algorithmic}
\Procedure{INIT2\_SOL}{$n,m,R$}
\State $R[1] \gets 1$
\State $R[n] \gets 1$
\State $sumE \gets 0$
\For{$i\gets 2, n-1$} 
   \State $R[i] \gets 0$
\EndFor
\State $e = m-n+1$
\State $i\gets n-2$
\While {$(i > 1) \wedge (e > 0)$}
    \If {$i - 1 <= e$}
        \State {$R[i] \gets 1$}
        \State {$e \gets e - i+1$}
        \State {$sumE \gets sumE + i-1$}
    \EndIf
    \State $i\gets i-2$
\EndWhile
\State $i\gets n-1$
\While {$(i > 1) \wedge (e > 0)$}
    \If {$i - 1 <= e$}
        \State {$R[i] \gets 1$}
        \State {$e \gets e - i+1$}
        \State {$sumE \gets sumE + i-1$}
    \EndIf
    \State $i\gets i-1$
\EndWhile
\If {$(sumE != e)$}
   \State $i\gets 1$
   \State $done \gets 0$
   \While {$(i<n-2) \wedge (done == 0)$}
      \If {$R[i] == 1 \wedge R[i+1] == 0$}
        \State {$R[i] \gets 0$}
        \State {$R[i+1] \gets 1$}
        \State {$done \gets 1$}
      \EndIf
      \State $i\gets i+1$
   \EndWhile
\EndIf
\EndProcedure
\end{algorithmic}
}
\end{algorithm}

The procedure described in Algorithm~\ref{INIT2-pscode} enables to obtain a solution where zeros and ones alternate in the binary sequence ensuring a larger number of possibilities to transform this solution into its neighbors with respect to some selected neighborhoods. Consecutively, it is more suitable for population based BCOi approach that does not perform systematic search in any neighborhood, but modifies solutions by a random number of randomly selected transformations.

\subsubsection{Solution transformations}

Four transformations yielding to a feasible solution can be defined on the proposed solution representation. 
\begin{enumerate}
\item Each combination $\{1...01...10...1\}$ can be replaced with $\{1...10...01...1\}$
\item Analogously, $\{1...10...01...1\}$ could be replaced with $\{1...01...10...1\}$
\item If $r_i=1$, $r_j=r_k=0$, and $i=j+k$, then it is possible to modify this solution in such a way that 
$r_i=0$, $r_j=r_k=1$
\item Analogously, if $r_i=0$, $r_j=r_k=1$, and $i=j+k$, then it is possible to modify this solution in such a way that
$r_i=1$, $r_j=r_k=0$
\end{enumerate}

The first two transformations preserve the number of zeros and ones in the resulting solutions, while the last two enable to modify (increase or decrease) these numbers, however, without violating their sum, i.e., the number of edges. To simplify the implementation, we grouped the first two transformations into a single neighborhood, and the remaining two into the second neighborhood. In such a way, GVNS will use two neighborhoods within the VND procedure, while BCOi will randomly choose one of transformation each time it attempts to modify a solution from the population.

The output of the search is the Adjacency matrix (reconstructed from the $\mathbf{R}$ sequence corresponding to the found graph with maximum radius), as well as the spectrum of that graph. The search should be performed for all connected threshold graphs, i.e., for each $n$, all graphs with $m=n-1$ up to $m=n(n-1)/2$ should be examined. In addition, as both metaheuristics are stochastic search methods, it is necessary to perform at least 100 repetitions (with different seed value) for each particular graph with fixed values for $n$ and $m$.

\section{Experimental evaluation}
\label{exper}

Here we present the results of applying the proposed BCOi method to the problem of finding threshold graphs with maximum spectral radius. The main idea is to generate hypothesis about the structure of extremal graphs with given $n$ and $m$ that could be theoretically proven by the SGT experts.

\subsection{Testing environment}

The implementation of GVNS is performed in C++ and executed on Intel Xeon E5-2620 v3, 2.40GHz, 32 GB RAM Under Linux 4.19.12 and compiled with GCC 4.8.3. 
BCOi is coded in Dev-C++ V5.11 and executed on Intel i5-6200U 230GHZ, RAM of 12GB DDR4 under Windows 10.

For the fair comparison of results obtained with the two developed algorithms, we need to define the parameter values, seed settings, number of repetitions, stopping criterion and other data relevant for experimental evaluation.
To be able to control the experiment and to replicate the results, we used a fixed set of values for seed in GVNS and BCOi. For the sake of simplicity, seed value in the $i$-th execution equals $n*i+m$. We hope that 30 repetitions is enough to have a general judgement about the performance of the proposed metaheuristic approaches. Stopping criterion is set to 1000 evaluations of the objective function value. Namely, the most time consuming part of our approaches is the calculation of spectrum for each examined graph. In every iteration of BCOi, the spectrum is calculated $NC$ times for each of $B$ solutions in the population. Predicting the number of objective function calculation in GVNS is not that straightforward, as the number of feasible neighbors of a given solution in each of the considered neighborhoods may vary substantially. Therefore, we limited the execution of both GVNS and BCOi by the same number of function evaluations, regardless the algorithmic step when they occur.

In the preliminary experimental evaluations, parameter values for both methods are determined intuitively.
GVNS has a single parameter $k_{max}$ and its value is set to $n/2$. It is important to note that we are not always able to perform the desired number of transformations and that the distance between the starting and resulting solution in the Shaking step may be smaller than expected. Only second neighborhood is used in Shaking and the Local Search strategy is Best Improvement. The parameters of
BCOi are set as follows: $B=5$, $NC=10$, $o=rnd (n/3, 2n/3)$. Here, $o$ determines the number of transformations of a single solution in each forward pass and it is determined randomly from the given interval. The same remark about the distance between the initial and final solution holds for the forward pass in BCOi.  

As the test examples we used small instances with 8 vertices and the number of edges taking value from the interval [12,23] (out of 28 possible edges in the completely connected graph $K_8$), as well as some medium-size instances with 30 and 50 vertices. The number of edges was selected randomly from the middle of valid interval. Having in mind that, for a given number of vertices, the number of connected (threshold) graphs with respect to the number of edges is a bell-shape function, 
it was reasonable to examine the part of the interval that corresponds to the largest number of graphs, i.e., the largest search space.


\subsection{Results of experimental evaluation}

The comparison results of the proposed algorithms on the small instances are presented in Table~\ref{resO}.

\begin{table}[!htb]
\caption{Comparison of the GVNS and BCOi results on small graph examples}
\label{resO}       
\begin{center}
\begin{tabular}{|l|r|r|r|r|r|r|r|}
\hline
Graph  &\multicolumn{1}{c|}{Init.sol.} & \multicolumn{3}{c|}{GVNS} & \multicolumn{3}{c|}{BCOi}\\
\hline
 & Obj.val. & \#bests & best obj. & av. obj. & \#bests & best obj. & av. obj.\\
\hline
$G_{8,12}$  & 3.78& 30& 3.85& 3.85& 30&  3.85& 3.85\\
$G_{8,15}$  & 4.37& 30& 4.52& 4.52& 30&  4.52& 4.52\\
$G_{8,19}$  & 5.19& 30& 5.33& 5.33& 30&  5.33& 5.33\\
$G_{8,21}$  & 5.62& 30& 5.77& 5.77& 30&  5.77& 5.77\\
$G_{8,23}$  & 6.10& 30& 6.10& 6.10& 30&  6.10& 6.10\\
\hline
\end{tabular}
\end{center}
\end{table}

Table~\ref{resO} is organized as follows. The names of instances (containing the numbers of vertices and edges) are given in the first column. The second column contains the objective function value (spectral radius) of the initial solution obtained by the procedure presented in Algorithm~\ref{INIT2-pscode}. The next three columns contain the results provided by GVNS, number of best (out of 30 repetitions) graphs, the objective function value corresponding to the best graph, and average (over 30 trials) value of the spectral radius, respectively. The same data related to the execution of BCOi is provided in the remaining three columns of Table~\ref{resO}. 

As can be seen in Table~\ref{resO}, both algorithms were able not only to improve the objective function value of the initial solution, but also to provide the optimal solutions in all 30 executions. The exhaustive search over all threshold graphs with 8 vertices and the corresponding number of edges was performed to prove this. 

The binary sequences corresponding to the obtained results are as follows.

$\mathbf{R}(G_{8,12})=\{1\; 0\; 1\; 1\; 0\; 0\; 0\; 1 \}$

$\mathbf{R}(G_{8,15})=\{1\; 1\; 0\; 1\; 1\; 0\; 0\; 1 \}$

$\mathbf{R}(G_{8,19})=\{1\; 1\; 1\; 0\; 1\; 1\; 0\; 1 \}$

$\mathbf{R}(G_{8,21})=\{1\; 0\; 1\; 1\; 1\; 1\; 0\; 1 \}$

$\mathbf{R}(G_{8,23})=\{1\; 1\; 1\; 1\; 1\; 0\; 1\; 1 \}$

The results provided by the preliminary experiments on graphs with 30 and 50 vertices did not demonstrate the same quality and stability. Therefore, we performed some modifications and parameter tuning with an aim to increase the performance of both algorithms. In the new setup for GVNS, both neighborhood types are applied with the probability 0.5 in the shaking procedure. Local Search strategy is changed to First Improvement, i.e., as soon as an improved solution is obtained the current Local Search iteration is interrupted and the one starts from this solution. BCOi is modified in such a way that initial population in each (except the first) iteration contains 2 best-so-far solutions (intended for additional improvements and ensuring intensification of the search) and 3 completely new randomly generated solutions (that are needed to enable diversification of the search process). In addition, Finally, the stopping criterion is increased to 2000.

\begin{table}[!htb]
\caption{Comparison of the results obtained by GVNS and BCOi for graphs with 30 and 50 vertices}
\label{res1}       
\begin{center}
\begin{tabular}{|l|r|r|r|r|r|r|r|}
\hline
Graph &\multicolumn{1}{c|}{Init.sol.} & \multicolumn{3}{c|}{GVNS} & \multicolumn{3}{c|}{BCOi}\\
\hline
 & Obj.val. & \#bests & best obj. & av. obj. & \#bests & best obj. & av. obj.\\
\hline
$G_{30,100}$  & 10.96& 30& 12.34& 12.34      &  13&  12.34& 12.10\\
$G_{30,220}$  & 18.12& 30& 20.03& 20.03& 30&  20.03& 20.03\\
$G_{30,300}$  & 22.16& {\textbf{30}}& 23.65& {\textbf{23.65}}& 26&  23.65& 23.64\\
$G_{30,400}$  & 27.04& 30& 27.58& 27.58& 30&  27.58& 27.58\\
\hline
$G_{50,100}$  & 10.38& 30& 10.87& 10.87& 30&  10.38& 10.38\\
$G_{50,300}$  & 19.58& {\textbf{30}}& {\textbf{22.89}}& {\textbf{22.89}}& 25&  22.50& 22.19\\
$G_{50,500}$  & 26.80& {\textbf{30}}& {\textbf{30.33}}& {\textbf{30.33}}& 1&  30.18& 30.08\\
$G_{50,1000}$ & 41.87& 30& 44.02& 44.02& 30&  44.02& 44.02\\
\hline
\end{tabular}
\end{center}
\end{table}

The obtained results on larger graphs are presented in Table~\ref{res1}. This table has the same structure as the previous one.  Comparing the results from Table~\ref{res1} we can conclude that GVNS exhibits stable performance as it obtains the same result in all 30 repetitions. It outperforms the best BCOi results in 3 out of 8 tested graph examples. It is also evident that for small and large enough number of edges both algorithms perform equally good, which is an indication that those are the easiest cases. When the number of edges approaches the middle of the examined interval, the search space becomes larger and systematic search performed by GVNS yields better results that the random perturbations among population members in BCOi.

The binary sequences that correspond to the best obtained solutions are as follows.

$\mathbf{R}(G_{30,100})=\{1 1 1 1 1 1 1 0 1 1 1 1 1 0 0 0 0 0 0 0 0 0 0 0 0 0 0 0 0 1 \}$

$\mathbf{R}(G_{30,220})=\{1 1 1 1 1 1 1 1 1 1 1 1 1 1 1 1 1 1 1 0 1 0 0 0 0 0 0 0 0 1 \}$

$\mathbf{R}(G_{30,300})=\{1 1 1 1 1 0 1 1 1 1 1 1 1 1 1 1 1 1 1 1 1 1 1 1 0 0 0 0 0 1 \}$

$\mathbf{R}(G_{30,400})=\{1 1 1 1 1 1 1 0 1 1 1 1 1 1 1 1 1 1 1 1 1 1 1 1 1 1 1 1 0 1 \}$

$\mathbf{R}(G_{50,100})=\{1 1 1 1 0 1 1 1 1 1 1 0 0 0 0 0 0 0 0 0 0 0 0 0 0 0 0 0 0 0 0 0 0 0 0 0 0 0 0 0 0 0 0 0 0 0 0 0 0 1\}$

$\mathbf{R}(G_{50,300})=\{1 1 0 1 1 1 1 1 1 1 1 1 1 1 1 1 1 1 1 1 1 1 1 0 0 0 0 0 0 0 0 0 0 0 0 0 0 0 0 0 0 0 0 0 0 0 0 0 0 1\}$

$\mathbf{R}(G_{50,500})=\{1 1 1 1 1 1 1 1 1 1 1 1 1 1 0 1 1 1 1 1 1 1 1 1 1 1 1 1 1 1 1 0 0 0 0 0 0 0 0 0 0 0 0 0 0 0 0 0 0 1\}$

$\mathbf{R}(G_{50,1000})=\{1 1 1 1 1 1 1 1 1 1 1 1 1 1 1 1 1 1 1 1 1 1 1 1 1 1 1 1 1 1 1 1 1 1 1 1 1 1 1 0 1 1 1 1 1 0 0 0 0 1\}$

Analyzing the structure of final solution, the hypothesis can be stated that the largest spectral radius have graphs whose binary sequences contains 1s concentrated at the beginning. This contrasts the greedy procedure for creating initial solution. 

\section{Conclusion}
\label{concl}

We implemented GVNS and BCOi as the incomplete search for graphs with maximum radius. Promising results for further development are obtained by preliminary experimental evaluation. The main challenges are a large complexity of objective function computation and complex data structures used to store the solutions of the considered problem. This makes difficult to perform search on large graphs. As the possible directions for future research we identify:
\begin{itemize}
\item Optimization of the developed code
\item The possible inclusion of new neighborhoods/transformations
\item The analysis of possibility to utilize memory and learning from previously visited solutions
\item A careful parameter tuning
\end{itemize}

\medskip

\noindent {\bf Acknowledgements.} This work has been partially supported by the Serbian Ministry of Science, Technological Development, and Innovations Agreement No. 451-03-66/2024-03/200029, as well as by the Science Fund of the Republic of Serbia, Grant \#6767, Lazy walk counts and spectral radius of graphs—LZWK.

\bibliographystyle{splncs04}
\bibliography{bco,sgt,vns}

\end{document}